\documentclass[a4paper]{amsart}
\usepackage{pdfsync}
\usepackage{amssymb}
\usepackage{color}
\usepackage[initials]{amsrefs}

\theoremstyle{plain}
\newtheorem{theorem}{Theorem}%[section]

\newtheorem{proposition}[theorem]{Proposition}
\newtheorem{lemma}[theorem]{Lemma}
\newtheorem*{theorem*}{Theorem}
\newtheorem*{lemma*}{Lemma}
\newtheorem*{proposition*}{Proposition}
\newtheorem*{corollary*}{Corollary}

\theoremstyle{definition}

\newtheorem*{definition*}{Definition}
\newtheorem*{example*}{Example}
\newtheorem*{remark*}{Remark}

\begin{document}
\title{Exponential Sums by Irrationality Exponent}

%\begin{comment}
\author{Byungchul Cha}
\address{Department of Mathematics Education,
  Hongik University, 
  94 Wausan-ro, Mapo-gu,
Seoul 04066, Korea.}
\email{cha@hongik.ac.kr}

\author{Dong Han Kim}
\address{Department of Mathematics Education,
Dongguk University - Seoul, 30 Pildong-ro 1-Gil, Jung-gu, Seoul 04620, Korea.}
\email{kim2010@dgu.ac.kr}
%\end{comment}


%\date{\today}


\thanks{B.~Cha was supported by the Hongik University new faculty research fund. D.~H.~Kim was supported by National Research Foundation of Korea (RS-2023-00245719).
}

\begin{abstract} 
  In this article, we give an asymptotic bound 
  for the exponential sum of the M\"obius function
  $\sum_{n \le x} \mu(n) e(\alpha n)$ for a fixed irrational number $\alpha\in\mathbb{R}$.
  This exponential sum was originally studied by Davenport
  and he obtained an asymptotic bound of $x(\log x)^{-A}$ for any $A\ge0$.
  Our bound depends on the irrationality exponent $\eta$ of $\alpha$.
  If $\eta \le 5/2$, we obtain a bound of $x^{4/5 + \varepsilon}$  %%d
  and, when $\eta \ge  5/2$, our bound is $x^{(2\eta-1)/2\eta + \varepsilon}$.  %%d
  This result extends a result of Murty and Sankaranarayanan, who obtained the same bound in the case $\eta = 2$.
\end{abstract}

\maketitle

\section{Introduction}

%Let $e(\alpha)$ denote $e(\alpha)= \exp(2\pi\alpha\sqrt{-1})$. 
%Using Vinogradov's earlier work \cite{V}, Davenport \cite{D} established the following bound for the exponential sum: for a fixed (irrational) $\alpha$,
%\begin{equation}\label{VD}
%\sum_{n \le x} \mu(n) e(\alpha n) \ll x (\log x)^{-A},
%\end{equation}
%for any $A \ge 0$, where the implied constant depends only on $A$.
%Under Riemann hypothesis, Baker and Harman \cite{BH} improved the bound in \eqref{VD} as $O_{\varepsilon}(x^{\frac34 + \varepsilon})$.
%
%In this paper, we improve the bound of the above exponential sum in \eqref{VD} according to the \emph{irrationality exponent} of $\alpha$ and 
%Recall, from \cite{Bugeaud} for example, that the irrationality exponent of $\alpha$ is defined to be the supremum of the set of real numbers $\tau$ for which the inequality
%\begin{equation}\label{IrrataionalityMeasure}
%0 < \left|\alpha -\frac pq \right| < \frac1{q^{\tau}}
%\end{equation}
%is satisfied by infinitely many integer pairs $(p,q)$ with $q > 0$.
%The irrationality exponent is greater than or equal to 2 for any irrational $\alpha$.  %%d
%Our main theorem is 

%{\color{magenta} 
Let $e(\alpha)= \exp(2\pi\alpha\sqrt{-1})$. 
Using Vinogradov's earlier work \cite{V}, Davenport \cite{D} established the following bound for the exponential sum of the M\"obius function: for a fixed irrational $\alpha$,
\begin{equation}\label{VD}
  \sum_{n \le x} \mu(n) e(\alpha n) \ll_A x (\log x)^{-A},
\end{equation}
for any $A \ge 0$, where the implied constant depends only on $A$.

It is possible to improve the above bound at least conditionally, by introducing an assumption on a zero-free region of Dirichlet $L$-functions.
Hajela and Smith \cite{HS87} showed that
\begin{equation}
  \sum_{n \le x} \mu(n) e(\alpha n) \ll_{\varepsilon} x^{(a + 2)/3 + \varepsilon},
  \label{eq:Hajela_Smith}
\end{equation}
under the assumption that
\begin{equation}
  L(s, \chi) \text{ has no zeros in the half plane } \Re(s) > a.%   {\color{blue}\text{ what is $\sigma$?}}
  \label{eq:zero_free_assumption}
\end{equation}
Improving \eqref{eq:Hajela_Smith}, 
Baker and Harman \cite{BH} proved, under the assumption \eqref{eq:zero_free_assumption} for all Dirichlet $L$-functions, that
\begin{equation}
  \sum_{n \le x} \mu(n) e(\alpha n) \ll_{\varepsilon} x^{b + \varepsilon},
  \label{eq:B_bound}
\end{equation}
where
\[
  b =
  \begin{cases}
    a + \tfrac14 & \tfrac12 \le a <\tfrac{11}{20}, \\
    \tfrac45 & \tfrac{11}{20} \le a <\tfrac35, \\
    \tfrac{a+1}2 & \tfrac35 \le a <1. 
  \end{cases}
\]
Note that, if we assume GRH for all Dirichlet $L$-functions, the bound of Hajela and Smith becomes $x^{\frac56 + \varepsilon}$, 
while that of Baker and Harman is $x^{\frac34 + \varepsilon}$.
More recently, assuming \eqref{eq:zero_free_assumption} in the range $a\in [1/2, 4/7]$,
Zhang \cite{Zha} improved Baker and Harman's bound in this range by proving \eqref{eq:B_bound}
where
\[
  b = \frac{8a - 7a^2}{4 - 2a}.
\]
For more results on this topic, see \cite{JL19} and the introduction therein.

In this paper, we take a different approach by considering the \emph{irrationality exponent} of $\alpha$
without any assumption on a zero-free region of $L$-functions.
Recall, from \cite{Bugeaud} for example, that the irrationality exponent of $\alpha$ is defined to be the supremum of the set of real numbers $\tau$ for which the inequality
\begin{equation}\label{IrrataionalityMeasure}
  0 < \left|\alpha -\frac pq \right| < \frac1{q^{\tau}}
\end{equation}
is satisfied by infinitely many integer pairs $(p,q)$ with $q > 0$.
The irrationality exponent is greater than or equal to 2 for any irrational $\alpha$.  %%d
Our main theorem is 
%}
\begin{theorem}\label{thm:estimate}
  Suppose that $\alpha$ is a fixed irrational number whose irrationality exponent $\eta$ is finite. Let $\varepsilon > 0 $.
  If $2\le \eta \le \frac 52$, then
  \[
    \sum_{n \le x} \mu(n)e(\alpha n) \ll_{\varepsilon} 
    x^{\left(\frac45 + \varepsilon\right)} .
  \]
  If $\eta \ge \frac 52$, then
  \[ 
    \sum_{n \le x} \mu(n)e(\alpha n)  \ll_{\varepsilon}
    x^{ \left(\frac{2\eta -1}{2\eta} + \varepsilon\right)}.
  \]
\end{theorem}

When $\eta = 2$, that is, when $\alpha$ is of \emph{irrational type 1}, 
Theorem 1 has been already proved in 2002 by 
Murty and Sankaranarayanan (see Corollary 4 in \cite{MS02}).

Note that the bounds for the exponential sum given by the other authors we quoted above
are all uniform, in the sense that they are independent of $\alpha$. 
Thus what we achieve in this paper is to improve Davenport's bound at the expense of making it dependent on the irrationality exponent of $\alpha$. 

A bound for the exponential sum for the M\"obius function
can be used to prove Sarnak's M\"obius disjointness conjecture for the simplest    %%d
type of dynamical system, namely the rotation of $\mathbb{R}/\mathbb{Z}$ by $\alpha$. 
To recall the statement of the conjecture \cite{Sar} in general, let $(X, T)$ be a topological dynamical system. For any $x\in X$ and $f\in C(X)$, the conjecture asserts that, if $(X, T)$ is \emph{determinisitic}, namely, of entropy zero, then
\begin{equation}\label{TheConjecture}
  \sum_{n\le N} \mu(n) f(T^n x) = o(N)
\end{equation}
as $N\to\infty$.
If $X = \mathbb{R}/\mathbb{Z}$ and if $T=T_{\alpha}$ is the rotation by $\alpha$, that is, $T_{\alpha}(x) := x + \alpha$, the verification of \eqref{TheConjecture} boils down to showing the orthogonality of $\mu(n)$ against $e(\alpha n)$. 
In fact, one can easily deduce from \eqref{VD} that the left-hand side of \eqref{TheConjecture} is the big-$O$ of $N(\log N)^{-A}$ (whose implicit constant depends on $x, f$ and $A$), hence $o(N)$. 
This will in turn yield a better convergence rate with respect to the conjecture \eqref{TheConjecture} of Sarnak.   %%d   the sentence is moved
We note here that the conjecture \eqref{TheConjecture} has been proven for several other kinds of systems $(X, T)$, such as Thue-Morse sequences \cite{MR}, nilsequences \cite{GT}, and horocycle flows \cite{BSZ}.

Our proof begins with a few preliminary lemmas in \S\ref{sec:proposition}, culminating in Proposition~\ref{prop:VI}. 
The results in \S2 are originally due to Vinogradov \cite{V} and Vaughan \cite{Vau}. 
Our exposition mostly follows \cite{IK} and \cite{GT2}. 
The proof of Theorem~\ref{thm:estimate} is deduced from Proposition~\ref{prop:VI}.

\section{Lemmas}\label{sec:proposition}

We begin with the, so-called, \emph{Vaughan's identity} \cite{Vau}, following the exposition in Chapter 13 of \cite{IK} and Section 4 of \cite{GT2}: for fixed $M$ and $N$,

\begin{equation}\label{Vaughan}
  \sum_{m\le x}\mu(m) e(\alpha m) = 
  -T_{\mathrm{I}} + T_{\mathrm{II}} + O(\max\{ M, N \}),
\end{equation}
where
\[
  T_{\mathrm{I}} := \underset{\substack{\ell k \le x \\ k \le MN}}{\sum\sum} e(\alpha \ell k)  \sum_{\substack{n \mid k \\ n \le N}} \mu(n) \mu(k/n) 
\]
and
\[
  T_{\mathrm{II}} :=
  \underset{\substack{ kn \le x \\ k  > M, n  > N }}{\sum\sum} \mu(n) e(\alpha k n) \sum_{\substack{m \mid k \\ m > M}} \mu(m) 
  .
\]
To prove \eqref{Vaughan}, write
\[
  \mu(m) = \underset{bc | m}{\sum\sum} \mu(b)\mu(c)
\]
and break up the summation into four regions: 
(i) $b\le M$ and $c \le N$
(ii) $b\le M$ and $c > N$
(iii) $b > M$ and $c \le N$
(iv) $b > M$ and $c > N$.
Then it is easy to see that $\sum_{(\mathrm{ii})} = \sum_{(\mathrm{iii})}=-\sum_{(\mathrm{i})}$, if $m>\max\{M, N\}$.
Thus,
\begin{multline*}
  \sum_{m \le x} \mu(m) e(\alpha m ) = 
  - \underset{\substack{\ell m n \le x \\ m \le M, n \le N }}{\sum\sum\sum} \mu(m) \mu(n) e(\alpha \ell m n)  
  \\
  +
  \underset{\substack{\ell m n \le x \\ m > M, n > N}}{\sum\sum\sum} \mu(m) \mu(n) e(\alpha \ell m n) + O(\max\{M,N\}).
\end{multline*}
From this, a simple change of variable changes the first summation in the right hand side above as $T_{\mathrm{I}}$ and the second as $T_{\mathrm{II}}$, establishing \eqref{Vaughan}.

To bound $T_{\mathrm{I}}$ and $T_{\mathrm{II}}$, we write
\[ 
  T_{\mathrm{I}} =  \underset{\substack{ \ell k \le x \\ k \le MN}}{\sum\sum} e(\alpha \ell k) \gamma (k, N),
\]
and
\[ 
  T_{\mathrm{II}} =  
  \underset{\substack{ kn \le x \\ k > M, n > N }}{\sum\sum} \mu(n) e(\alpha kn) \tau(k,M),
\]
where 
\[
  \gamma (k, N):= \sum_{\substack{n \mid k \\ n \le N}} \mu(n) \mu(k/n), \qquad
  \tau(k,M) := \sum_{\substack{m \mid k \\ m > M}} \mu(m).
\]
Let $d(k):=\sum_{u|k} 1$ be the divisor function, that is, the number of positive divisors of $k$. Then it is well known that
\[
  d(k) \ll_\varepsilon k^{\varepsilon}.
\]
Moreover, we clearly have
\[
  | \gamma (k, N) | \le d (k), \qquad | \tau(k,N) | \le d (k).
\]

Now we appeal to the following two lemmas from \cite{IK}. Suppose that $q$ satisfies
\[
  \left|
  \alpha - \frac aq
  \right|
  \le
  \frac 1{q^2} \text{ with } (a, q) = 1.
\]
\begin{lemma}[Lemma 13.7 in \cite{IK}]\label{Lemma1}
  \[
    \sum_{1\le m \le M}
    \left|
    \sum_{1\le n \le \frac xm}
    e(\alpha mn)
    \right|
    \ll
    \left(
      M + \frac xq + q
    \right)
    \log(2qx)
  \]
  as $x\to\infty$.
\end{lemma}
This proves that 
\[
  T_{\mathrm{I}}  \ll_{\varepsilon}
  \left(
    MN + \frac xq + q
  \right)
  x^{\varepsilon} \log(2qx).
\]
For $T_{\mathrm{II}}$, we use 
\begin{lemma}[Lemma 13.8 in \cite{IK}]\label{Lemma2}
  Let $\{\alpha_n\}$ and $\{\beta_m\}$ be sequences of complex numbers whose absolute values are bounded by 1. Then,
  \[
    \underset{\substack{mn \le x \\ m > M, n>N}}{\sum\sum} 
    \alpha_m \beta_n e(\alpha mn)  
    \ll
    \left(
      \frac xM + \frac xN + \frac xq + q
    \right)^{\frac12}
    x^{\frac12} (\log x)^2,
  \]
  as $x\to\infty$.
\end{lemma}
Likewise,
\[
  T_{\mathrm{II}}  \ll_{\varepsilon}
  \left(
    \frac xM + \frac xN + \frac xq + q
  \right)^{\frac12}
  x^{\frac12 + \varepsilon} (\log x)^2.
\]
We summarize what we proved as 
\begin{proposition}\label{prop:VI}
  Let $\varepsilon>0$ be given. Then, for fixed $M$ and $N$
  \begin{equation}%\label{Vaughan}
    \sum_{m\le x}\mu(m) e(\alpha m) = 
    -T_{\mathrm{I}} + T_{\mathrm{II}} + O(\max\{ M, N \}),
  \end{equation}
  where $T_{\mathrm{I}}$ and $T_{\mathrm{II}}$ satisfy
  \begin{equation}\label{TypeI}
    T_{\mathrm{I}}  \ll_{\varepsilon}
    \left(
      MN + \frac xq + q
    \right)
    x^{\varepsilon} \log(2qx)
  \end{equation}
  \begin{equation}\label{TypeII}
    T_{\mathrm{II}}  \ll_{\varepsilon}
    \left(
      \frac xM + \frac xN + \frac xq + q
    \right)^{\frac12}
    x^{\frac12 + \varepsilon} (\log x)^2.
  \end{equation}
\end{proposition}
\section{Proof of Theorem~\ref{thm:estimate}}

Let $\eta$ be the irrationality exponent of $\alpha$ and fix a real $\tau $ with $\tau > \eta$. Then by the definition of irrationality exponent, 
\[
  \frac1{q^{\tau}} <
  \left|
  \alpha - \frac pq
  \right|
\]
is satisfied for all but finitely many integer pairs of $p$ and $q$. Then there exists $i_0 = i_0(\tau, \alpha)$ such that, for all $i\ge i_0$,
\[
  \frac1{q_{i-1}^{\tau}} <
  \left|
  \alpha - \frac{p_{i-1}}{q_{i-1}}
  \right|
  < \frac1{q_iq_{i-1}},
\]
where $p_i/q_i$ is the $i$-th principal convergent in the continued fraction of $\alpha$. This implies
\begin{equation}\label{Qgrowth}
  q_i < {q_{i-1}}^{\tau - 1}
\end{equation}
for all $i \ge i_0$. Now, given large $x$, choose $i$ so that 
${q_{i-1}}^{\tau - 1} \le x^{(\tau-1)/\tau} < {q_i}^{\tau-1}$. Then $q:=q_i$ satisfies
\begin{equation}\label{Xrange}
  q^{\frac{\tau}{\tau - 1}} < x < q^{\tau}.
\end{equation}
From this, we deduce
\begin{equation}\label{XQrange}
  q <  x^{1 - \frac1{\tau}}, \qquad \text{and} \qquad  \frac xq  < x^{1 - \frac1{\tau}}.   %%d
\end{equation}
We now invoke Proposition~\ref{prop:VI} with $M=N=x^{\frac25}$. From \eqref{XQrange} and \eqref{TypeI},
\begin{equation}\label{TI}
  T_{\mathrm{I}} \ll_{\varepsilon}
  \left(
    x^{\frac 45} + x^{1-\frac1{\tau}}
  \right)
  x^{\varepsilon},
\end{equation}
and also from \eqref{TypeII}
\begin{equation}\label{TII}
  T_{\mathrm{II}} \ll_{\varepsilon}
  \left(
    x^{\frac 85} + x^{2-\frac1{\tau}}
  \right)^{\frac12}
  x^{\varepsilon}.
\end{equation}
Combining \eqref{TI} and \eqref{TII}, we complete the proof of Theorem~\ref{thm:estimate}.

\section*{Acknowledgements}
The authors are thankful to Peter Humphries who brought \cite{MS02} to the authors' attention.

%\bibitem{BeNa} V. Berth\'e and H. Nakada, \emph{On continued fraction expansions in positive characteristic: equivalence relations and some metric properties}, Expo. Math. \textbf{18} (2000), 257--284.

\end{document}